\documentclass{amsart}

\newtheorem{theorem}{Theorem}[section]

\newtheorem{lemma}[theorem]{Lemma}

\newcommand{\BB}{{\mathbb B}}
\newcommand{\ZZ}{{\mathbb Z}}
\newcommand{\RR}{{\mathbb R}}
\newcommand{\CC}{{\mathbb C}}

\newcommand{\PP}{{\mathbb P}}
\newcommand{\HH}{{\mathbb H}}
  
\newcommand{\II}{{\mathbb I}}

\newcommand{\cP}{{\mathcal P}}

\newcommand{\cJ}{{\mathcal J}}

\newcommand{\la}{\lambda}

\newcommand{\diag}{\rm{diag}}

\title[ Perron-Frobenius-Ruelle theorem and two weight Hilbert transform
]{  Noncommutative Perron-Frobenius-Ruelle theorem, two weight 
Hilbert transform, and almost periodicity  }
\author{A. Volberg, P. Yuditskii}
  \thanks{Partially supported by NSF grant DMS-0200713,  
  the grant for IAS and 
  the Austrian Science
Found FWF, project number: P16390--N04}

\begin{document}

\maketitle

\begin{abstract}
We consider the Jacobi matrix generated by a balanced measure of
hyperbolic polynomial map. The conjecture of Bellissard says that this
matrix should have an extremely strong periodicity property. We show
how this conjecture is related to a certain noncommutative version of
Bowen--Ruelle theory, and how the two weight Hilbert transform
naturally appears in this context.
\end{abstract}

\section{Introduction and Main results}
\label{Intr}

Let
$f$
be an expanding polynomial with real Julia set $J(f)$, $\deg f= N$.
We recall that $J(f)$ is a nonempty compact set of points which do not go to infinity under forward iterations of $f$.
Under the normalization 
$$
f^{-1}:[-\xi,\xi]\to [-\xi,\xi]
$$
such a polynomial is well defined by position
of its critical values 
$$
\{t_i=f(c_i): f'(c_i)=0,\  c_i>c_j \ \text{for}\ i>j
\}.
$$
Expanding, or hyperbolic polynomials are those, for which 
$$
c_i\notin J(f), \forall i\,.
$$
The term ``expanding" is deserved because for expanding polynomials one has the following inequality
\begin{equation}
\label{expanding}
\exists Q >1,\,\,|(f^n)'(x)| \geq c Q^n, \forall x\in J(f)\,.
\end{equation}
Here and in everything that follows $f^n$ means $n$-th iteration of $f$, $f^n = f\circ f\circ....f$.

We will always use letter $T$ for polynomial $f^n$, $\deg T =N^n$. We will always use letter $d$ for this degree, $d=N^n$.

We will say that $f$ is sufficiently hyperbolic if 
$$
\forall i\,\,\,{\rm dist} (f(c_i), J(f)) \geq A\,,
$$
with a sufficiently large $A$.

Let us mention that for $f$ with a {\it real}
Julia set one has $|f(c_i)|>\xi$ since all solutions
of $f(x)=\pm \xi$ should be real.

We recall now Perron-Frobenius-Ruelle (PFR) theorem in a form convenient for us.
Let $\phi$ be a $Hol(\alpha)$ function on $J(f)$. We define the Perron-Frobenius-Ruelle (PFR) operator
$$
\mathcal{L}_{\phi} = \mathcal{L}_{\phi,f}: C(J(f)) \rightarrow C(J(f))
$$ as follows
$$
\mathcal{L}_{\phi} \psi(x) :=\sum_{\lambda: f(\lambda) =x} e^{\phi(\lambda)}\psi(\lambda)\,.
$$

PFR theorem states that if $\rho$ denotes the spectral radius of this operator then
$$
\rho^{-n} \mathcal{L}_{\phi}^n \psi (x) \rightarrow h(x)\int \psi(y) d\nu(y)\,,
$$
where $h$ is the unique eigenvector of $\mathcal{L}_{\phi}$ with eigenvalue $\rho$, $\nu$ is is the unique eigenvector of $\mathcal{L}_{\phi}^*$ with eigenvalue $\rho$. Moreover, $h$ is H\"older continuous if $\phi$ is H\"older continuous.

Let us emphasize that the requirement on $f$ is just  to be expanding (hyperbolic).

We reformulate this result now. Or, rather we will formulate its essential part in a different form.
In fact, it turns out that at the heart of this result lies the following theorem. See \cite{Bo}, for example.

\begin{theorem}
\label{0.1}
Let $f$ be hyperbolic. Let $\phi \in Hol_{\alpha}(J(f)),  \psi \in Hol_{\alpha}(J(f)),$ then there exist $C<\infty, q\in (0,1), \gamma>0$ (all depending only on $\alpha$) such that
\begin{equation}
\label{0.2}
\Bigl|\frac{\mathcal{L}_{\phi}^n \psi (x_1) }{\mathcal{L}_{\phi}^n 1(x_1) } -\frac{\mathcal{L}_{\phi}^n \psi (x_2) }{\mathcal{L}_{\phi}^n 1 (x_2) }\Bigr| \leq C\,q^n|x_1-x_2|^{\gamma}\,.
\end{equation}
\end{theorem}

Consider the operator $G_{\phi}= G_{\phi,f}$ acting on probability measures  on $J(f)$ by the formula
$$
G_{\phi}\mu = \frac{\mathcal{L}_{\phi}^*\mu}{\|\mathcal{L}_{\phi}^*\mu\|} =\frac{\mathcal{L}_{\phi}^*\mu}{\langle 1,\mathcal{L}_{\phi}^*\mu\rangle}\,.
$$
Then \eqref{0.2} means

\begin{equation}
\label{0.3}
|\langle \psi , G_{\phi}^n\delta_{x_1} \rangle -\langle \psi , G_{\phi}^n\delta_{x_2} \rangle| \leq  C\,q^n|x_1-x_2|^{\gamma}\,.
\end{equation}
This is for any test function $\psi \in Hol_{\alpha}(J(f))$, under assumption that $\phi \in Hol_{\alpha}(J(f))$.

In what follows $\phi$ will have the following form
$$
\phi := -t \log|f'|\,, t\in \mathbb{R}\,.
$$

\noindent{\it Example}. When $t=0$ we have that $G^n \delta_x$ is a sum
of delta measures with charges $1/d=1/N^n$ located at all $T$-preimages
of $x$.
We would like to understand the PFR theorem as a consequence of a certain fact of noncommutative nature. Sometimes this is indeed so, we can prove this. We would like to prove this noncommutative
fact always, for all hyperbolic dynamics, but we can manage only  a weaker estimate.

To explain what noncommutative proposition we have in mind, let us notice that there is a natural operator for which $G_{\phi}^n\delta_{x}$ is a spectral measure. This is just Jacobi matrix built by this probability measure. Let us recall that to built the Jacobi matrix by a probability measure $d\mu(\lambda)$  with support on the real line, one just orthogonalizes polynomials with respect to this measure, and Jacobi matrix is the matrix of multiplication by independent variable $\lambda$ written in the basis of orthonormal polynomials.  

So let $T=f^n$,  $\deg T=d =N^n$, $J(x)=J_T(x)$ be a Jacobi matrix built by measure
$$
\mu_x= G_{\phi, f}^n\delta_x = G_{\phi, T} \delta_x\,,
$$
where $\phi=-t \log|f'|$ (as always in what follows).

We already explained that $J(x)$ is canonically defined. 
Another way to define $J(x):\mathbb{C}^d\rightarrow \mathbb{C}^d$ 
is to write
$$
\langle (z-J(x)^{-1} e_0,e_0\rangle = \int \frac{d\mu_x(\lambda)}{z-\lambda}=
\sigma \sum_{k=1}^d \frac{e^{\phi(\lambda_k)}}{z-\lambda_k}\,,
$$
where $\lambda_1(x),....,\lambda_d(x)$ are all $T$-preimages of $x$ .

We can consider $J(x)$ as a result of applying PFR operator to $1\times 1$ matrix $x$.
We want to prove the following claim, which deserves to be called {\it a noncommutative Perron-Frobenius-Ruelle (PFR) theorem}.
\begin{equation}
\label{0.4}
\|J(x_1) -J(x_2)\| \leq C\,q^n |x_1-x_2|^{\gamma}\,.
\end{equation}

We cannot prove \eqref{0.4} for all hyperbolic $f$ and all $\phi$. But here are our main results.

\begin{theorem}
\label{lip}
Let $f$ be hyperbolic.  Let $x_1, x_2\in J(f)$.  Let $T=f^n$, we build $J(x)=J_T(x)$ using $\phi=-t \log|f'|$ , $0\leq t\leq2$. Then there exists $C$ such that independently of $n$
\begin{equation}
\label{lipeq}
\|J(x_1) -J(x_2)\| \leq C |x_1-x_2|\,.
\end{equation}
\end{theorem}

\begin{theorem}
\label{contraction}
Let $f$ be sufficiently hyperbolic.   Let $x_1, x_2\in J(f)$. Let  $T=f^n$, we build $J(x)=J_T(x)$ using $\phi=0$. Then there exists $c<1$ such that independently of $n$
\begin{equation}
\label{contractioneq}
\|J(x_1) -J(x_2)\| \leq c^n |x_1-x_2|\,.
\end{equation}
\end{theorem}

\vspace{.2in}

\noindent{\bf Acknowledgements.} The authors are grateful to Michael Shapiro,  Fedja Nazarov and Sergei Treil for valuable discussions. The first author acknowledge with deep gratitude the grant from IAS that allowed to use the
stimulating atmosphere of this institution.

\section{PFR theorem from noncommutative PFR theorem}
\label{PFRfromPFR}

This uniform closeness of operators quite easily imply the type of weak closeness of Theorem \ref{0.1}, \eqref{0.2} or \eqref{0.3}. We just use the fact that measures in, say, \eqref{0.3}, are spectral measures of $J(x_i)$, and we use the Jackson-Bernstein type theorem to approximate the H\"older continuous function $\psi$ by functions holopmorphic in narrowing neighborhoods of $J(f)$ with the speed
$\varepsilon^{\tau}$, where $\varepsilon$ is the width of a neighborhood, and $\tau$ depends on $\alpha$ in H\"older property of $\psi$. This is very standard, and it shows that really PFR theorem can be considered as a consequence of a noncommutative claim---\eqref{0.4}.

\section{ Faybusovich--Gehktman flow}
\label{GF}
 
Let $J(x)=J_T(x), T=f^n$, and we think now that $n$ is large but 
fixed, and we are heading towards the proof of Theorem \ref{lip} with
constant independent of $n$. 
If we think that $x$ is the time, then the flow of Jacobi matrices $J(x)$ of size $d\times d$ ($d=N^n$)
can be treated {\it alike} the Toda flow. But {\it unlike} the 
Toda flow, the spectrum $sp(J(x))$ is not time
independent, that is, it is not $x$-independent. This brings a modification to the equation of Toda flow.
Such modifications were considered by  Faybusovich and Gehktman in
\cite{GF}, we grateful to M. Shapiro who indicated this to us.  
Let $()^{\cdot}$ denote the differentiation with respect to ``time" $x$. We write $J$ instead of $J(x)$ for brevity. Recall that in the standard basis $e_0, e_1,....,e_{d-1}$ of $\mathbb{C}^d$ the matrix of $J$ is three diagonal. Recall that $sp(J)$ is always equal  $\lambda_1(x),....,\lambda_d(x)$  (all $T$-preimages of $x$). If $g$ is a continuous function on this spectrum, we know what is $g(J)$ by functional calculus of self-adjoint operators. We need one more definition.
Given $x$ we consider the orthonormal polynomials $P_0(\lambda;x) =1, P_1(\lambda;x),..., P_{d-1}(\lambda;x) $ of degrees $0,1,..., d-1$ correspondingly. They are orthonormal with respect to measure
$$
d\mu_x(\lambda) := \sigma\sum_{k=1}^d e^{\phi(\lambda_k(x)} \delta_{\lambda_k(x)}\,,
$$
where from now on always $\phi =-t\log|f'|$, $1/\sigma = \sum_{k=1}^d e^{\phi(\lambda_k(x)} $.

\begin{theorem}
\label{flow}
Our $J(x)$ satisfies the nonlinear ODE
\begin{equation}
\label{floweq}
J^{\cdot} = F(J) + [G,J]\,,
\end{equation}
where $F(J)$ is a function of a self-adjoint operator $J=J(x)$ with function $F:= (T')^{-1}$ on $sp(J)$.
Operator $G$ is skew self-adjoint, and its matrix in the standard basis has upper triangular part $G_-$ equal to
\begin{equation}
\label{G}
G_- = (DF(J) + \frac12\phi'(J)F(J))_-\,.
\end{equation}
Finally operator $D$ is given by the formula
\begin{equation}
\label{D}
\langle D^* e_k,e_m\rangle := \sigma \sum_{k=1}^d e^{\phi(\lambda_k(x)} P_k'(\lambda_i(x)) P_m(\lambda_i(x))\,.
\end{equation}
\end{theorem}

The form of $D$ will allow to prove
\begin{theorem}
\label{D}
\begin{equation}
\label{Deq}
[J,D] = I - c\cdot \langle \cdot, F^{-1} e_0\rangle e_{d-1}\,.
\end{equation}
And therefore
\begin{equation}
\label{DFeq}
[J, DF]e = Fe , \forall e\,\text{orthogonal}\, \text{to}\, e_0\,.
\end{equation}
\end{theorem}

In Toda flow $F=0$ and $G_-= R(J)_-$ for a certain function $R$. Here
we are in a more complicated situation, but let us observe
\begin{equation}
\label{smallF}
\|F(J)\| \leq Cq^n, \,\,q <1\,.
\end{equation}
In fact, just use \eqref{expanding}. Then if $x\in J(f)$ we have $sp(J)\subset J(f)$ and \eqref{expanding} implies automatically the latter inequality \eqref{smallF}. If $x$ is not on $J(f)$, but is separated from the critical values of $T$, inequality \eqref{expanding} also holds on $T$-preimages of $x$, and this set is exactly $sp(J)$.

It is very good that $F(J)$ is small as in \eqref{smallF} because to prove Theorem \ref{lip} it 
is  hence enough to prove

\begin{theorem}
\label{commutator}
\begin{equation}
\label{commutatoreq}
\|[G,J]\| \leq C\,,
\end{equation}
where $C$ is independent of $n$.
Then automatically
\begin{equation}
\label{deriveq}
\|J^{\cdot}\|\leq C'\,.
\end{equation}
\end{theorem}

\noindent{\bf Remark.} We cannot prove that $G$ is uniformly bounded, moreover this looks to be false.

\section{Uniform boundedness  of the commutator}
\label{unif}

We postpone the explanation of Theorem \ref{flow}. Now we will take it for granted, and we prove Theorem \ref{commutator}, which of course immediately gives the first main result---Theorem \ref{lip}.

Let 
$$
H:= DF(J) + \frac12\phi'(J) F(J)\,.
$$
Let us adopt the following convention. If we write
$$
A=B + small\,,
$$ 
we always mean that the small term is at most $C q^n$, with $q<1$ in norm.

Then, for every $x$
$$
[G,J] = \text{three diaginal}  + small\,,
$$
just see \eqref{floweq}, \eqref{smallF}.

This means that to prove the uniform boundedness of the commutator is the same as to prove the uniform boundedness of its matrix elements in the standard basis. Moreover,
it is enough to prove the uniform boundedness of ``three diagonal" elements only:
$$
|\langle[G,J]e_m, e_{m-1}\rangle|, |\langle[G,J]e_{m}, e_{m}\rangle|, |\langle[G,J]e_{m-1}, e_{m}\rangle|\,.
$$
Moreover the skew symmetry of $G$ implies that $[G,J]$ is self-adjoint, so only the first and the second 
type of elements should be checked.

Recall that the upper triangular parts coincide:
$$
G_- =H_-\,.
$$
Let us denote the diagonal elements of $J$ by $a_0=a_0(x),..., a_{d-1} = a_{d-1}(x)$,
and let below diagonal elements be $b_1=b_1(x),..., b_{d-1} = b_{d-1}(x)$.

Let us write the following equalities ($G_+$ is the lower triangular part of skew symmetric $G$)
$$
\langle[G,J]e_m, e_{m-1}\rangle =\langle[G_-,J]e_m, e_{m-1}\rangle+ \langle[G_+,J]e_m, e_{m-1}\rangle
=\langle[G_-,J]e_m, e_{m-1}\rangle\,.
$$
In fact using only the lower triangular property of $G_+$ we see that $ G_+ J e_m\in span(e_m,...)$ and so is orthogonal to $e_{m-1}$. Also,
$ G_+  e_m\in span(e_{m+1},...)$ and so is orthogonal to $Je_{m-1}$.

Now we conclude that
$$
\langle[G,J]e_m, e_{m-1}\rangle =\langle[H_-,J]e_m, e_{m-1}\rangle=
\langle[H,J]e_m, e_{m-1}\rangle-
$$
$$
\langle[H_0,J]e_m, e_{m-1}\rangle-\langle[H_+,J]e_m, e_{m-1}\rangle\,,
$$
where $H_+$ is the lower triangular part of $H=DF +\phi'(J)F$, and $H_0$ is its diagonal part
(that may exist because nobody said that $H$ is skew symmetric). We
claim that this means
\begin{equation}
\label{H1}
\langle[G,J]e_m, e_{m-1}\rangle =-\langle[H_0,J]e_m, e_{m-1}\rangle +small\,.
\end{equation}
In fact,
the last term $\langle[H_+,J]e_m, e_{m-1}\rangle$ equals zero--we checked this for any lower triangular matrix.
On the other hand, the term $\langle[H,J]e_m, e_{m-1}\rangle$ is equal to 
$\langle[DF +\phi'(J)F(J),J]e_m, e_{m-1}\rangle=\langle[DF ,J]e_m, e_{m-1}\rangle=
-\langle Fe_m, e_{m-1}\rangle$. The last equality follows from \eqref{DFeq}. So this term is small in our sense, and \eqref{H1} is proved.

We can now use
$$
Je_m = b_m e_{m-1} + a_m e_m + b_{m+1} e_{m+1}\,
$$
to plug it into \eqref{H1} and get it rewritten 
\begin{equation}
\label{H2}
\langle[G,J]e_m, e_{m-1}\rangle =b_m\langle He_m,e_m\rangle - b_{m-1} \langle He_{m-1},e_{m-1}\rangle +small\,.
\end{equation}

We can also calculate $\langle[G,J]e_m, e_{m}\rangle$.
\begin{equation}
\label{H3}
\langle[G,J]e_m, e_{m}\rangle =2b_{m+1}\langle He_m,e_{m+1}\rangle -2 b_{m} \langle He_{m-1},e_{m}\rangle +small\,.
\end{equation}

In fact, skew symmetry of $G$ implies $\langle[G,J]e_m, e_{m}\rangle  = 2 \langle[G_-,J]e_m, e_{m}\rangle  = 2\langle[H_-,J]e_m, e_{m}\rangle $.  Again replace $H_-$ by $H-H_0 -H_+$.
Then
$$
\langle[G,J]e_m, e_{m}\rangle = 2\langle[H,J]e_m, e_{m}\rangle-2\langle[H_0,J]e_m, e_{m}\rangle-
2\langle[H_+,J]e_m, e_{m}\rangle=: 2A-2B-2C\,.
$$
Obviously $B=0$.  To see that $A=small$ let us use \eqref{Deq}: $[J,D]Fe =Fe-c\cdot \langle Fe, F^{-1}e_0\rangle e_{d-1}$. In other words $[J,DF] e_m = Fe_m -c\cdot\langle e_m, e_0\rangle_{d-1}$.
But we saw that $[H,J] =[DF,J]$. Therefore $A= \langle[H,J]e_m,e_m\rangle = small $ for any $m=0,1,..., d-1$. We are left to see what is $C$.
$$
C= \langle H_+ J e_m, e_m\rangle - \langle H_+ e_m, Je_m\rangle=
$$
$$
\langle H_+ ( b_m e_{m-1} + a_m e_m + b_{m+1} e_{m+1}), e_m\rangle -
\langle H_+ e_m, b_m e_{m-1} + a_m e_m + b_{m+1} e_{m+1}\rangle=
$$
$$
b_m \langle H_+ e_{m-1}, e_m\rangle - b_{m+1}\langle H_+ e_m, e_{m+1}\rangle\,.
$$
In both expressions here one can replace $H_+$ by $H$ without changing these scalar
products. So \eqref{H3} is proved.

From \eqref{H2} and \eqref{H3} it follows that the estimate  of the norm of the commutator $[G,J]$
follows from the estimate of operator $H=DF(J) + \frac12\phi'(J) F(J)$.

If we can prove the uniform boundedness of $H$, we can prove, therefore, Theorem \ref{commutator} and, thus Theorem \ref{lip}.

\section{The uniform boundedness of $H=DF(J) + \frac12\phi'(J) F(J)$. Two weight Hilbert transform}
\label{twoweight}

To prove the uniform boundedness of $H$ we need to understand $D$ better. To do this we will write $H$ in a different basis, and we will see that $DF$  becomes a two weight Hilbert transform (almost).
Then we use our knowledge of the boundedness of two weight Hilbert transform. This will prove the uniform boundedness of $H$, and, as a result, will prove Theorem \ref{commutator} and Theorem \ref{lip}.

We already introduced polynomials orthonormal with respect to $\mu_x$.
Now consider the following matrices:

\begin{equation}
\cP=
\begin{bmatrix}
P_0(\la_1)&...&P_0(\la_d) \\
  \vdots&&\vdots & \\
 P_{d-1}(\la_1)&...&P_{d-1}(\la_d) 
\end{bmatrix}.
\end{equation}

\begin{equation}
\Phi=
\begin{bmatrix}
e^{\phi(\la_1)}&&& \\
 & \ddots\\
&&& e^{\phi(\la_d)}
\end{bmatrix}.
\end{equation}

Recall that $1/\sigma=\sum_{k=1}^d e^{\phi(\la_k)}$.
Put
\begin{equation}
\label{PP}
\PP := \sqrt{\sigma}\cdot\cP\cdot\Phi\,.
\end{equation}
The orthonormality of polynomials $\{P_k\}$ with respect to $d\mu_x =\sigma\sum_{k=1}^d e^{\phi(\la_k)}\delta_{\la_k}$  means that matrix $\PP$ is an orthogonal matrix.

Let $\Lambda$ be 
$$
\Lambda:= \cP^{-1}\, J \,\cP, \,\,\, R:= \cP^{-1}\, DF\,\cP\,.
$$
Then it easy to see that also
$$
\Lambda:= \PP^{-1}\, J \,\PP\,.
$$
This just because $\cP$ and $\PP$ are different only up to a diagonal matrix.
Moreover,
$$
F(\Lambda) = \PP^{-1}\,F(J) \,\PP =\cP^{-1}\, F(J) \,\cP\,.
$$

\begin{lemma}
\label{LaR}
Then
\begin{equation}
\Lambda=
\begin{bmatrix}
\la_1&&& \\
 & \ddots\\
&&& \la_d
\end{bmatrix}.
\end{equation}
\begin{equation}
\label{LaReq}
[\Lambda, R]=
\begin{bmatrix}
0,\frac1{T'(\la_1)},...,\frac1{T'(\la_1)}\\
\frac1{T'(\la_2)}, 0,..., \frac1{T'(\la_2)}\\
\vdots\\
\frac1{T'(\la_d)},...,\frac1{T'(\la_d)}, 0
\end{bmatrix}
\end{equation}
\end{lemma}

\begin{proof}
The first formula is obvious as the unitary (it is even orthogonal) matrix $\PP$  transforms $J$ into its diagonal form $\Lambda$. To see the second formula we have to notice that 
\begin{equation}
\label{LaR1}
[\Lambda, R]=\cP^{-1}[J,DF]\cP = F(\Lambda)  - c\cdot \langle \cdot, \cP^*e_0\rangle \cP^{-1} e_{d-1}\,.
\end{equation}
This is because of formula \eqref{Deq}.

Now it is obvious by definition that 
$$ \cP^* e_0 =
\begin{bmatrix}
1\\ \vdots\\1\end{bmatrix}
$$
And let us see that 
\begin{equation}
\label{LaR2}
c\cdot \cP^{-1} e_{d-1} = 
\begin{bmatrix}
\frac 1{T'(\la_1)}\\ 
\vdots\\ \frac1{T'(\la_d)}\end{bmatrix}
\end{equation}
This and \eqref{LaR1} will finish the lemma.
To prove \eqref{LaR2} let us notice that denoting
$$
c\cdot \cP^{-1} e_{d-1} = 
\begin{bmatrix}
v_1\\ \vdots\\ v_d\end{bmatrix}
$$
we obtain from \eqref{LaR1} and the form of $\cP^* e_0 $ that
$$
[\Lambda, R] = F(\Lambda) -
\begin{bmatrix}
v_1,...,v_1\\
\vdots\\
v_d,...,v_d
\end{bmatrix}
$$
But $\Lambda$ is diagonal, and so the diagonal of the LHS vanishes. But this gives us $v_i=$
$i$-th diagonal element of $F(\Lambda)$, which is $\frac1{T'(\la_i)}$. So \eqref{LaR2} is proved, and
lemma is shown.

\end{proof}

Recall that
$$R:= \cP^{-1}\, DF\,\cP\,.
$$
\begin{lemma}
\label{R}
The matrix elements of $R$ are as follows
$$
r_{ij} = \frac1{T'(\la_i)} \frac1{\la_i-\la_j},\,\,\text{if}\,\, i\neq j\,.
$$
$$
r_{ii} = \frac12\frac{T''(\la_i)}{(T^{'}(\la_i))^2}\,.
$$
\end{lemma}

\begin{proof}
The non-diagonal terms can be immediately read from \eqref{LaReq} of the previous lemma.
On the other hand
$$
R^* \begin{bmatrix}1\\ \vdots\\1\end{bmatrix} =0\,.
$$
In fact, 
$$
R^* \begin{bmatrix}1\\ \vdots\\1\end{bmatrix} =\cP^* FD^*(\cP^*)^{-1}\begin{bmatrix}1\\ \vdots\\1\end{bmatrix}\,,
$$
But we know that
$$
(\cP^*)^{-1}\begin{bmatrix}1\\ \vdots\\1\end{bmatrix} =e_0,
$$
and $D^*e_0=0$ (from \eqref{D}). So our first equality is proved.
It means that the sum of column elements of $R$ is zero for 
every column. We knew all elements of $R$ except the diagonal ones, but
this sum property gives us the diagonal elements too. An easy residue
theorem application gives the formula $r_{ii} =
\frac12\frac{T''(\la_i)}{(T^{'}(\la_i))^2}\,.$
\end{proof}

Let 
$$ 
K:= \PP^{-1} DF \PP\,.
$$
Compare this with $R=\cP^{-1} DF \cP\,.$  As $\cP$ and $\PP$ are almost
the same matrices---the  difference is in diagonal factor we get
\begin{theorem}
\label{K}
\begin{equation}
\label{Keq}
K= \diag\left(\frac{e^{-\phi(\la_i)/2}}{T'(\la_i)}\right)
\begin{bmatrix}
 \frac12\frac{T"(\la_1)}{T^{'}(\la_i)},...,\frac1{\la_1-\la_d}\\
 \vdots\\
\frac1{\la_1-\la_d},..., \frac12\frac{T"(\la_i)}{T^{'}(\la_d)}\end{bmatrix}diag(e^{\phi(\la_i)/2})
\end{equation}
And matrix $H$ is unitary equivalent to
\begin{equation}
\label{Heq}
\diag\left(\frac{e^{-\phi(\la_i)/2}}{T'(\la_i)}\right)
\begin{bmatrix}
 \frac12(\log|T^{'}|+\phi)'(\la_i),...,\frac1{\la_1-\la_d}\\
 \hdots\\
\frac1{\la_1-\la_d},..., \frac12(\log|T'|+\phi)'(\la_d)\end{bmatrix}\diag(e^{\phi(\la_i)/2})
\end{equation}
In particular, having in mind that $\phi =-t\log|T'|$ we get
that 
matrix $H$ is unitary equivalent to
\begin{equation}
\label{Heqt}
\HH=\HH_t=\diag\left(\frac{1}{|T'(\la_i)|^{1-\frac{t}2}}\right)
\begin{bmatrix}
 \frac{1-t}2(\log|T^{'}|)'(\la_i),...,\frac1{\la_1-\la_d}\\
 \hdots\\
\frac1{\la_1-\la_d},..., \frac{1-t}2(\log|T'|)'(\la_d)\end{bmatrix}
\diag\left(\frac1{T'(\la_i)|^{\frac{t}2}}\right)
\end{equation}
\end{theorem}

\begin{proof}
The first relation \eqref{Keq} follows from Lemma \ref{R} and from formula \eqref{PP} that relates
$\cP$ and $\PP$ via a multiplication on a  diagonal  matrix.
But $\PP^{-1} H\PP = K + \frac12\PP^{-1} \phi'(J)F(J)\PP= K +\frac12 F(\Lambda) \phi'(\Lambda)$,
and \eqref{Heq} follows.
\end{proof}

We are ready to prove the uniform boundedness of $H$.

\begin{theorem}
\label{HH}
Let $f$ be a hyperbolic polynomial of degree $N$ with real Julia set $J(f)$. Let $T=f^n$, $\deg T=d=N^n$.
Let $\{\la_1,...,\la_d\}$ be all $T$-preimages of $x\in J(f)$. 
Matrix $\HH_t$ is uniformly bounded independently of $n$, $x$, and $t, 0\leq t\leq 2$.
Therefore, so is matrix $H$.
\end{theorem}

We already saw that the proof of this theorem finishes the proof of Theorem \ref{commutator}, and thus, of our first main result, Theorem \ref{lip}.

\begin{proof}
The diagonal part $\HH^0 := \frac{1-t}{2}\diag (\log|T'|)'(\la_i))$ is bounded uniformly in $n$ and $x\in J(f)$ just by Koebe distortion theorem, it is a standard fact depending only on hyperbolicity of $f$.
(Notice that for $t=1$ this matrix vanishes!)
Let us consider now the ``out-of-diagonal" part 
$$
\BB=\BB_t = \HH_t - \HH_t^0\,.
$$
Consider the counting measure on $\{\la_1,...,\la_d\}$: 
$d n=d n_x =\sum_{k=1}^d \delta_{\la_k}$.
Now we can notice easily that $\BB^*$ is unitary equivalent to the following integral operator
$$
g\in L^2(d n) \rightarrow |T'(x)|^{-\frac{t}2}\int\frac1{x-y} \frac{|T'(y)|^{\frac{t}2}}{|T'(y)|}\,g(y)\,d n(y)\in L^2(d n)\,.
$$
Changing variable $f:= g\cdot |T'|^{1-\frac{t}2}$ and changing measure $d \nu(y) := \frac{d n(y)}{|T'(y)|^{2-t}}$ we come to a unitary equivalent operator
$$
f\in L^2(d \nu) \rightarrow |T'(x)|^{-\frac{t}2}\int\frac1{x-y} \,f(y)\,d \nu(y)\in L^2(d n)\,.
$$
Put
$$
d\kappa := |T'(x)|^{-t}\,d n(x)\,.
$$
The norm of $\BB$ is equal to the norm of the two weight Hilbert transform
$$
H_{\nu} f := \int_{y\neq x} \frac1{x-y} \,f(y)\,d \nu(y): L^2(\nu) \rightarrow L^2(\kappa)\,.
$$

The story of two weighted problems in Harmonic analysis is beyond the scope of this work. We will just choose the result  convenient for our narrow purpose here. A paper \cite{TVZ} looks like  specially written for the occasion. However, the reader who wants to familiarize her/himself with two weighted estimates is referred to \cite{Vo}  and to the vast literature cited there. We just make two remarks. First one is that the two weight estimates for operators with positive kernels is more or less well understood due to the works of Eric Sawyer (many of them are cited in \cite{Vo}). On the other hand the singular kernel two weight estimates are not completely understood even for the simplest singular kernels (like the Hilbert transform). There is only one kernel--the dyadic singular kernel corresponding to the Martingale transform, where the technique of Bellman function gives a full criterion of boundedness.
See \cite{NTV2}. There is no ``classical" approach to this so far. And if kernel becomes just slightly more complicated than the dyadic one (for example the Hilbert transform) there is no real understanding.
(The criterion of Cotlar-Sadosky \cite{CS}  is very nice but its language seems to be not applicable here.)
Some criterion which ``seems to be" the right one is considered in the last two chapters of 
\cite{Vo}. There are some counterexamples to other ``right criteria" in \cite{NV}. 

But we have to find a certain applicable and easily verifiable sufficient condition of two weight boundedness of the Hilbert transform. The question is very intimately related to a so-called problem of Sarason: describe when the product of two Toeplitz operators is bounded. Dechao Zheng found a wonderful sufficient condition in \cite{DZ}. It was then adopted in \cite{TVZ} to two weight Hilbert transform. One of the main results of \cite{TVZ} will be applied here---it is perfect for our goals.

\vspace{.1in}

Let us introduce notations. The symbol $\langle f\rangle_I$ will denote the usual averaging $\frac1{|I|}\int_I f\,dx$, where
$I$ is an interval on a real line.
The symbol $P_I f$ denotes the Poisson averaging, namely, $\frac1{\pi}\int_{\RR} \frac{|I|}{(x-c)^2 + |I|^2} f(x)\,dx$, where $c$ is the center of $I$. In other words it is the value of the Poisson extension of $f$ at the point $c+i\cdot |I|\in \CC_+$.

We prove Theorem \ref{HH} if we prove the following result.
\end{proof}

\begin{theorem}
\label{twHt}
The norm of 
$$
H_{\nu} : L^2(\nu) \rightarrow L^2(\kappa)
$$
is uniformly bounded in $n$, $x$, and $t, 0\leq t\leq 2$.
\end{theorem}
 
\begin{proof}
To prove it we need the following result
\begin{lemma}
\label{gauge}
Let $udx, vdx$ be two positive measures on the line. Let  $g(t) =|t|^{1+\varepsilon}$, with $\varepsilon>0$.
If for every interval $I$ we have
\begin{equation}
\label{HMW}
P_I g(u) \cdot P_I g(v) \leq C
\end{equation}
with $C<\infty$ independent of $I$, then the two weight Hilbert transform
$$
H_{udx} : L^2(udx) \rightarrow L^2(vdx)
$$
is bounded, and its norm depends only on $C<\infty$ and $\varepsilon>0$.
\end{lemma}

\noindent{\bf Remark.} The reader may wonder why we need the gauge function $g$ here?
It turns out that $P_I u \cdot P_I v\leq C$ is not sufficient for the boundedness of the Hilbert transform in general. See \cite{N}, or \cite{NV}.

\vspace{.1in}

Let us reduce Theorem \ref{twHt} to this lemma. 

We will to this in two stages. Our first goal will be to prove the following weaker version of \eqref{HMW}:
\begin{equation}
\label{HMWw}
\langle g(u)\rangle_I\langle g(v)\rangle_I \leq C\,.
\end{equation}

\vspace{.2in}

Let us replace 
$$
d\,\nu(y)= \frac{d n(y)}{|T'(y)|^{2-t}}
$$ by $u(y)dy$, where
\begin{equation}
\label{u}
u(y):= \sum_{i=1}^d \frac1{|T'(y)|^{1-t}}\chi_{I_i}(y)\,,
\end{equation}
where $I_i$ is the $i$-th preimage of $[-\xi,\xi]$ under $T$ (left to right).

Similarly replace
$$
d\,\kappa(x)= \frac{d n(x)}{|T'(x)|^{t}}
$$ by $vdx$, where
\begin{equation}
\label{v}
v(x):= \sum_{i=1}^d \frac1{|T'(x)|^{t-1}}\chi_{I_i}(x)\,,
\end{equation}

\begin{lemma}
\label{antidiscr}
The norm of $H_{udx} : L^2(udx) \rightarrow L^2(vdx)$ bounds the norm of
$H_{d\nu} : L^2(d\nu) \rightarrow L^2(d\kappa)$.
\end{lemma}
\begin{proof}
Intervals $I_i$ are separated as their ``centers":
$$
dist(I_i, I_j) \asymp |\la_i-\la_j|\,.
$$
The constants of equivalence depend only on hyperbolicity. Every test function $f\in L^2(\nu)$, $f=(f_1,...,f_d)$ can be replaced by $F:= \sum f_i\chi_i$, and clearly
$$
\|f\|_{L^2(\nu)} \asymp \|F\|_{L^2(udx)}
$$
as
$$
|I_i| \asymp \frac1{|T'(\la_i)|}\,.
$$
\end{proof}
Now we prove the following 
\begin{lemma}
\label{DZ}
Let $g(t) =|t|^{1+\varepsilon}$.  Let $u,v$ be as in \eqref{u}, \eqref{v}. Then
$$
\sup_I \langle g(u)\rangle_I\langle g(v)\rangle_I<\infty\,.
$$
\end{lemma}
\begin{proof}
Let us think from now on that $\xi=1$. Let us first consider the case of the ``largest" interval: $I=[-\xi,\xi]=[-1,1]$. (We can consider only this or smaller intervals as the supports of all measure in question are in side this interval.)
Then
$$
\langle udx\rangle_I\langle vdx\rangle_I \asymp (\sum_{k-1}^d \frac1{|T'(\la_k)|^{2-t}})(\sum_{k-1}^d \frac1{|T'(\la_k)|^{t}})\,.
$$
Notice that
\begin{equation}
\label{Green}
 \frac1{|T'(\la_k)|} \leq\frac{C}{d}\,,
 \end{equation}
 where $d$ (as always) is $\deg T, T=f^n$.
 In fact, the LHS is equivalent to the distance to $I_k$ of the $k$-th component  of $T^{-1}(\Gamma)$, where $\Gamma$ is the circle of radius, say, $2\xi=2$ centered at zero.  The constants of equivalence depend only on the hyperbolicity of $f$. This is just Koebe distortion theorem. The RHS is equivalent to the value on the $k$-th component   of $T^{-1}(\Gamma)$ of Green's function of $\Omega = \CC\setminus T^{-1}([-3/2,3/2])$. In fact, this Green's function is
 $$
 \frac1d\log\Bigl|\frac23 T(z) + \sqrt{\frac49 T^2(z) -1}\Bigr| \asymp \frac1d,\,\text{if}\,\, z\in T^{-1}(\Gamma)\,.
 $$
 On the other hand Green's functions grow of the domain grows. So
 $$
 G_{\Omega}(z) \geq G_{\CC_+}(z) \asymp \Im z \geq c\cdot dist(z, I_k)\geq c\cdot  \frac1{|T'(\la_k)|}\,,
 $$
 if $z$ is on the top of $k$-th component  of $T^{-1}(\Gamma)$.
 
 Now \eqref{Green} and the Cauchy inequality give ($0\leq t<2$)
 $$
 \left(\sum_{k=1}^d \frac1{|T'(\la_k)|^{2-t}}\right)
 \left(\sum_{k=1}^d
\frac1{|T'(\la_k)|^{t}}\right)\leq
  d^{1-t}d^{1-\frac{2-t}{2}} \left(\sum_{k=1}^d \frac1{|T'(\la_k)|^2}
\right)^{1-\frac{t}2}\leq
 $$
\begin{equation}
\label{sumsum1}\left
 (d\cdot \sum_{k=1}^d \frac1{|T'(\la_k)|^2} \right)
 ^{1-\frac{t}2}\,
 \end{equation}
 For $t=2$
 \begin{equation}
 \label{sumsum2}
 \left(\sum_{k=1}^d \frac1{|T'(\la_k)|^{2-t}}\right)
 \left(\sum_{k=1}^d \frac1{|T'(\la_k)|^{t}}\right)\leq
 d\cdot \sum_{k=1}^d \frac1{|T'(\la_k)|^2} \,.
 \end{equation}
 \begin{lemma}
 \label{pressure2}
 \begin{equation}
 \label{pressure2eq}
 d\cdot \sum_{k=1}^d \frac1{|T'(\la_k)|^2} \leq C d^{-\tau}\,,
 \end{equation}
 where $C<\infty, \tau>0$ depend only on hyperbolicity of $f$.
 \end{lemma}
 \begin{proof}
 Again we use \eqref{Green} to get
 $$
 d\cdot \sum_{k=1}^d \frac1{|T'(\la_k)|^2} \leq C\cdot \sum_{k=1}^d \frac1{|T'(\la_k)|}\,.
 $$
 The last expression is equivalent to the length of $T^{-1}([-1,1])$ (and hence is bounded independently of $d=N^n$). To see our better estimate \eqref{pressure2eq} we shall recall the notion of pressure.
 For hyperbolic dynamics $f$ one introduces {\it the pressure}
 \begin{equation}
 \label{pressure}
 P(t) := \lim_{n\rightarrow\infty}\frac1n\log\Bigl(\sum_{k=1}^{N^n} \frac1{|(f^n)'|^t(\la_k)}\Bigr)\,.
 \end{equation}
 Here $\la_k, k=1,..., N^n$ are all preimages of a point. The limit exists, and gives us a convex and strictly decreasing function on $-\infty<t<\infty$.
 If we have a convention that $\log$ is in base $N$, we also have
 $$
 P(0)=1\,.
 $$
 It is known that the only root of $P$ is $\delta=Hdim J(f)$. As the dynamic is hyperbolic $\delta<1$.
 So
 $$
 P(1) =-\tau<0\,.
 $$
 This proves the lemma.
 \end{proof}
 The result of the lemma can be written as follows
 \begin{equation}
 \label{P02}
 P(0) + P(2) < 0\,.
 \end{equation}
 Our elementary inequalities \eqref{sumsum1} and \eqref{sumsum2} show that a more general fact is true:
 \begin{equation}
 \label{between0and2}
 P(t) + P(2-t) <0,\,\, 0\leq t \leq 2\,.
 \end{equation}
 Actually, it is trivial to see that \eqref{between0and2} follows from \eqref{P02} for any convex function $P$.
 Now notice that continuity of the pressure implies
 \begin{equation}
 \label{between0and2better}
 P(t) + P(2-t) <0,\,\, -\epsilon\leq t \leq 2+\epsilon\,.
 \end{equation}
 This is for a small positive $\epsilon$. We will need this now very much.
 
 Let us again consider the case of the ``largest" interval: $I=[-\xi,\xi]=[-1,1]$. 
 (We can consider only this or smaller intervals as the supports of all
measure in question are in side this interval.) But now we will use the
gauge function $g(t)= |t|^{1+\epsilon}$ with precisely this
$\epsilon$--the one from \eqref{between0and2better}.
$$
\langle u^{1+\epsilon}dx\rangle_I
\langle v^{1+\epsilon}dx\rangle_I \asymp \left(\sum_{k=1}^d
\frac1{|T'(\la_k)|^{2-t'}}\right)
\left(\sum_{k=1}^d \frac1{|T'(\la_k)|^{t'}}\right)\,.
$$
Here $t' = t+t\epsilon -\epsilon$. Of course, the range of $t'$ is $[-\epsilon, 2+\epsilon]$ and \eqref{between0and2better} gives
\begin{equation}
\label{theMain}
\langle u^{1+\epsilon}dx\rangle_I\langle v^{1+\epsilon}dx\rangle_I \leq C\cdot d^{-\eta},\,\,\text{for some positive}\,\,\eta\,.
\end{equation}
Here the interval was the largest possible. What if we take subintervals of $[-1,1]$?

First let us consider only ``dyadic" intervals $I$. We call the interval the interval $I$ dyadic if there exists
$m, 0\leq m\leq n$ such that $I= $ a component of $(f^m)^{-1}([-1,1])$. Such intervals form the set  $D_m$
of dyadic intervals of rank $m$. $D=\cup_{m=0}^n D_m$. We call them ``dyadic" even though they are not.
Recall that we ``smeared" our measures $d\nu$, $ d\kappa$ over intervals of $D_n$.
So let us fix $m, 0\leq m\leq n$ and an interval $I\in D_m$.  Let us split $T=f^n$ as follows: $T= f^{n-m}\circ f^m=: T_2\circ T_1$ on $I\in D_m$. Then 
$$
T_1(I) =[-1,1]\,,
$$
and 
\begin{equation}
\label{235}
|I|^{-1} \asymp |T_1'(x)|,\,\, x\in I\,.
\end{equation}
We want to estimate
$$
\langle u^{1+\epsilon}dx
\rangle_I\langle v^{1+\epsilon}dx\rangle_I \asymp\frac1{|I|}
\left(\sum_{i:\la_i\in I}\frac1{|T'(\la_i)|^{2-t'}}\right)
\frac1{|I|}
\left(\sum_{i:\la_i\in I} \frac1{|T'(\la_k)|^{t'}}\right)\,.
$$
Let $\{\mu_j\}_{j=1}^{N^{n-m}}$ be $f^{n-m}$ preimages of $x$.
Notice that
$$
\forall i: \la_i \in I\,\,\exists j\,\,T_1\la_i =\mu_j\,.
$$
Call $d_2 = N^{n-m}$. The expression we want 
to estimate is, by chain rule, and \eqref{235}  bounded by
$$
c\cdot \left(\sum_{j=1}^{d_2} \frac
1{|T_2'(\mu_j)|^{2-t'}}\right)\cdot
\left(\sum_{j=1}^{d_2}
\frac1{|T_2'(\mu_j)|^{t'}}\right).
$$
But this is bounded by $C\cdot d_2^{-\eta}$ by 
\eqref{between0and2better}. This is obtained by exactly the same
reasoning as we get \eqref{theMain}. Only $d_2$ replaces $d$. 
So far we proved \eqref{theMain} only for all ``dyadic" intervals. 
Let $J_n(f):= \cup_{i\in D_n} I$. If  for {\it any} interval
$I_0\subset [-1,1]$ such that $I_0\cap J_n(f) \neq \emptyset$  we would
have that there exists a ``dyadic" interval of comparable length that
contains $I_0\cap J_n(f)$, then \eqref{theMain} for $I=I_0$ would
follow from \eqref{theMain}  for ``dyadic" intervals. For usual dyadic
intervals this is of course false. It is obvious that one cannot always
find the dyadic interval of comparable length containing a given
interval. But in our situation this is true. 
\begin{lemma}
\label{dyadic}
Let $I_0\subset [-1,1]$ such that $I_0\cap J_n(f) \neq \emptyset$. 
Let $I$ denote the smallest interval from $D$ containing $I_0\cap
J_n(f)$. Then
\begin{equation}
\label{242}
|I_0|\geq c\cdot |I|\,,
\end{equation}
where $c>0$ depends only on hyperbolicity of $f$.
\end{lemma}
\begin{proof}
Along with dyadic intervals $D_1$ we have the collection of gap intervals $G_1$ between them.
Preimages of intervals of $G_1$ give gaps $G_k, k=2,...,n$. Take our $I_0$. Let $J$ be a gap interval inside it. If there is none then  $I_0\cap J_n(f)$ coincides with one interval of $D_n$.
 And \eqref{242} holds.  So let $J\subset I_0$ of the smallest generation $k=1,...,n$. $J\in G_k$.
 Then it lies in a dyadic interval $I\in D_{k-1}$. 
 Let us prove that
 \begin{equation}
\label{245}
I_0\cap J_n(f) \subset I\,,
\end{equation}
Interval $I$ has one or two neighbors of generation $k-1$ or smaller generation $m< k-1$. If \eqref{245} is false then $I_0$ should intersect one of these neighbors. But then it should contain the gap of generation $\ell\leq k-1$. This contradicts the choice of $J$. So \eqref{245} holds. 

One of the branches of $f^{-(k-1)}$ maps $[-1,1]$ onto $I$. Call this branch $g_{k-1}$. Moreover, univalently
$$
g_{k-1} : U\rightarrow U_I\,,
$$
where U is an open topological disc containing $[-1,1]$, and $I\subset U_I$.
Also $g_{k-1}$ maps a gap $L\in G_1$ onto $J$. Now Koebe distortion theorem implies
\begin{equation}
\label{248}
\frac{|J|}{|I|} \geq c\cdot \frac{|L|}{|[-1,1]|}\geq c_1 >0\,.
\end{equation}
Here $c,c_1$ depend only on $f$, not on $n$.
Obviously \eqref{248} implies \eqref{242}. Lemma \ref{dyadic} is proved.

\end{proof}

Together with \eqref{theMain} for dyadic intervals (already shown) it gives \eqref{theMain} for all intervals.
This proves \eqref{HMWw}. This is almost the proof of Lemma \ref{DZ}.

\vspace{.1in}

But to finish the proof of this lemma we need to pass from \eqref{HMWw}, which we have just proved to \eqref{HMW}. 

\vspace{.1in}

To do that we need still a couple of lemmas. First notations. Let 
$$
t'=t + \epsilon t -\epsilon, \, 0\leq t\leq 2\,.
$$
$$
\tau_0 := \tau_0 (t,\epsilon)=-[P(t') + P(2-t')] \,.
$$

\begin{lemma}
\label{eqviv}
Let ``dyadic" interval $I$ belong to $D_k$. Then
\begin{equation}
\label{eqviveq}
\langle u^{1+\epsilon}dx\rangle_I
\langle v^{1+\epsilon}dx\rangle_I \asymp N^{ -\tau_0 (n-k)}\,.
\end{equation}
\end{lemma}

\begin{proof}

\begin{equation}
\label{eqviveq1}
\langle u^{1+\epsilon}dx\rangle_I
\langle v^{1+\epsilon}dx\rangle_I \asymp \left(\sum_{j=1}^{N^{n-k}}
\frac1{|T'(\la_j)|^{2-t'}}\right)
\left(\sum_{j=1}^{N^{n-k}} \frac1{|T'(\la_j)|^{t'}}\right)\,.
\end{equation}

On the other hand, the hyperbolicity of dynamics $f$ standardly implies more than the existence of the limit in \eqref{pressure}. The more is actually known. Namely,
\begin{equation}
\label{pressureeqviv}
\sum_{\la: f^n \la =x} \frac1{|(f^n)'|^t(\la)}\asymp e^{-P(t) n}\,,
\end{equation}
where the constants of comparison do not depend on $n$ or $x\in J(f)$.

Now \eqref{eqviveq1} and \eqref{pressureeqviv} give us \eqref{eqviveq}, and the lemma is proved.

\end{proof}

\begin{lemma}
\label{doubling}
Let $I_m\in D_m$ be a ``dyadic" interval.  Let $I_{m+1}$ be its ``dyadic
subinterval of  $D_{m+1}$. Then
\begin{equation}
\label{doublingeq}
\int_{I_m} u^{1+\epsilon}\,dx \geq (1+\delta) \int_{I_{m+1}} u^{1+\epsilon}\,dx\,,
\end{equation}
where $\delta>0$ is independent of $m$.
\end{lemma}

\begin{proof}
Let us denote $I_{m+1}$ by $L$, and let $K$ be another interval from 
$D_{m+1}$ lying inside $I_m$.
Then using repeatedly the estimates from below and from above in Lemma \ref{eqviv} we get
$$
\int_{I_m} u^{1+\epsilon}\,dx \leq C\cdot N^{-\tau_0(n-m)} |I_m|\frac1{\langle v^{1+\epsilon}\rangle_{I_m}}\leq
$$
$$
C\cdot N^{-\tau_0(n-m)} |I_m|^2\frac1{\int_{I_m} v^{1+\epsilon}}\leq C\cdot N^{-\tau_0(n-m)} |I_m|^2\frac1{\int_{K} v^{1+\epsilon}}\leq
$$
$$
C\cdot N^{-\tau_0(n-m)} |I_m| \frac{|I_m|}{|K|}
\frac{\langle u^{1+\epsilon}\rangle_K}
{\langle u^{1+\epsilon}\rangle_K\langle v^{1+\epsilon}\rangle_K}
\leq
$$
$$
C\cdot \frac{|I_m|^2}{|K|^2}\int_K u^{1+\epsilon}\,dx \leq C\int_Ku^{1+\epsilon}\,dx\,.
$$
Of course we used in the last line that the lengths of an interval of $D_m$ and its ``son" from $D_{m+1}$ are comparable.
Let us rewrite the previous inequality as follows
$$
\int_Ku^{1+\epsilon}\,dx \geq \frac1C \int_{I_m}u^{1+\epsilon}\,dx \,.
$$
Then 
$$
\int_{I_m}u^{1+\epsilon}\,dx \geq \frac1C \int_{I_m}u^{1+\epsilon}\,dx + \int_{L}u^{1+\epsilon}\,dx\,.
$$
Lemma \ref{doubling} follows.
\end{proof}

\begin{lemma}
\label{PI}
Let $I\in D_k$ then
\begin{equation}
\label{PIeq}
P_I v^{1+\epsilon} \leq C\cdot N^{-\tau_0(n-k)} \frac1{\langle u^{1+\epsilon}\rangle_I}\,.
\end{equation}
\end{lemma}

\begin{proof}

Let us denote $I:=I_0\subset I_1\subset I_2 \subset.... \subset I_k =[-1,1]$ the nest of ``dyadic" intervals so that $I_j\in D_{k-j}$, $j=0,1,...,k$.
Then
$$
P_I v^{1+\epsilon} \leq C\sum_{j=0}^{k} \frac{|I|}{|I_j|} \langle v^{1+\epsilon}\rangle_{I_j}\,.
$$
Using Lemma \ref{eqviv} we can continue
$$
P_I v^{1+\epsilon} \leq C\sum_{j=0}^{k} \frac{|I|}{|I_j|} \frac{N^{-\tau_0(n-(k-j))}}{\langle u^{1+\epsilon}\rangle_{I_j}}\,.
$$
Now using Lemma \ref{doubling} we rewrite this 
$$
P_I v^{1+\epsilon} \leq C\sum_{j=0}^{k} |I|\frac{N^{-\tau_0(n-(k-j))}}{(1+\delta)^j\int_I u^{1+\epsilon}\,dx}\,.
$$
And this is exactly \eqref{PIeq}. Lemma is proved.
  \end{proof}
  
To finish the proof of Lemma \ref{DZ} we apply the previous lemma to $v$ and to $u$ and write
$$
P_I u^{1+\epsilon} P_I v^{1+\epsilon}  \leq C \frac{N^{-2\tau_0(n-k)} }{\langle u^{1+\epsilon}\rangle_I\langle v^{1+\epsilon}\rangle_I}\,.
$$  
We are left to use the estimate from below part of Lemma \ref{eqviv} to get

$$
P_I u^{1+\epsilon} P_I v^{1+\epsilon}  \leq C \cdot N^{-2\tau_0(n-k)},\,\,\forall\,I\in D_k\,.
$$
In Lemma \ref{DZ} one needs this same estimate but for {\it every subinterval} $I$  of $[-1,1]$.
Fix such an interval.  Consider first the case: there exists $\ell \in J_n(f)$ such that  $I\cap \ell \neq \emptyset$.  
Such an $I$ can be ``big" or small".  Big means $C\cdot I $ contains an interval of $\ell$.  Otherwise, $I$ is ``small" and then it intersects only one $\ell\in D_n$ and
$$
|II \leq |\ell|\,.
$$
In the latter case, we use the fact that $u$ is constant on $\ell$ to write
$$
P_I u^{1+\epsilon} \leq C\langle u^{1+\epsilon} \rangle_I  + P_{\ell}u^{1+\epsilon} \,.
$$
And the same for $v$.

But using again the fact that  $u$ is constant on $\ell$ to write
$$
P_I u^{1+\epsilon} \leq A\langle u^{1+\epsilon} \rangle_{\ell}  + B P_{\ell}u^{1+\epsilon} \leq CP_{\ell}u^{1+\epsilon}\,.
$$
And the same for $v$.

So in this case $P_I u^{1+\epsilon} P_I v^{1+\epsilon} \leq CP_{\ell} u^{1+\epsilon} P_{\ell} v^{1+\epsilon} $, and this has been proved to be universally bounded (depending only on hyperbolicity) for $\ell\in D$.

Now suppose that $I$ is ``big". 
Let $J$ be a ``dyadic" interval (that is $J\in D:= \cup_{k=0}^{n}D_k$) of maximal length contained in $C\cdot I$.  It is easy to see from Lemma \ref{dyadic} that
$$
|J|
\geq a\cdot |I|\,\ a>0\,.
$$
Here $a$ is independent of $I$. Then Harnack inequality and this previous relationship show that
\begin{equation}
\label{intersect}
P_I u^{1+\epsilon} P_I v^{1+\epsilon}  \leq C \cdot N^{-2\tau_0(n-k)},\,\,\forall\,I\in D_k
\end{equation}
for such an $I$ too. 

Finally, if $I\cap J_n(f) = \emptyset$ we can find the smallest interval $L$ such that $L\cap J_n(f)\neq\emptyset$ and $I\subset L$. Notice that $u=v=0$ on $L$ except the endpoint(s).
Therefore,
$$
P_I u^{1+\epsilon}  \leq C P_L u^{1+\epsilon} ,\, P_I v^{1+\epsilon}  \leq C P_L v^{1+\epsilon}\,.
$$
But for the interval $L$ \eqref{intersect} has just been proved.
Therefore it holds for $I$ too.

\vspace{.1in}

Lemma \ref{DZ} is completely proved.
\end{proof}

In its turn, it prove Theorem \ref{twHt}.

\end{proof}

Our first main result--Theorem \ref{lip}--is completely proved. 

\vspace{.2in}

Our second main result--Theorem \ref{contraction}--can be found in \cite{PVY1} or below.

\section{Sufficiently large hyperbolicity and contractivity of noncommutative PFR map for $\phi=0$. Almost periodicity.}
\label{contractivity}

Here we will discuss Theorem \ref{contraction} proved in \cite{PVY1}. Moreover, we will explain that not only
$$
\|J(x_1) - J(x_2)\|\leq c^n\,|x_1-x_2|\,, \,\,c<1\,,
$$
but that also there is an operator analog of this fact
$$
\|\cJ(J_1) - \cJ(J_2)\|\leq c^n\,|J_1-J_2|\,,\,\,c<1\,.
$$

But there will be restrictions. First of all $\phi=0$ only (so $t=0$ only). Secondly, unlike the previous sections, where it was not very important that $f$ is a polynomial, here this will be very much used. And thirdly at last, not just hyperbolicity, but only sufficiently large hyperbolicity allows us to prove this contractivity.

We do not know whether it is true in general. Or for other $t$'s, $\phi$'s.

A Jacobi matrix $J:l^2(\ZZ)\to l^2(\ZZ)$ is called almost periodic if
the family
$$
\{S^{-k}JS^k\}_{k\in \ZZ},
$$
where $S$ is the shift operator, $S|k\rangle=|k+1\rangle$,
is a precompact in the operator topology.

\begin{proof}[Example]
Let $G$ be a compact abelian group, $p(\alpha), q(\alpha)$ be continuous
functions on $G$,  $p(\alpha)\ge 0$. Then $J(\alpha)$ with the
coefficient sequences $\{p(\alpha+k\mu)\}_k, \{q(\alpha+k\mu)\}_k$,
$\mu\in G$, is almost periodic.
\end{proof}

Let us show that in fact this is a general form of almost periodic
Jacobi matrices. For a given almost periodic $J$ define the metric on
$\ZZ$ by
$$
\rho_J(k):=||S^{-k}JS^k-J||.
$$
Evidently $\rho_J(k+m)\le\rho_J(k)+\rho_J(m)$.  Then $J=J(0)$,
where $G=I_J$, $I_J$ is the closer
of $\ZZ$ with respect to $\rho_J$, and $\mu=1\in I_J$.

Recall that for a given system of integers $\{d_k\}_{k=1}^\infty$
one can define
 the set 
 of $d$ adic numbers as
\begin{equation}
\II=\underleftarrow{\lim}\{\ZZ/d_1...d_k\ZZ\},
\end{equation}
 that is $\alpha\in \II$ means that $\alpha$ is
 a sequence $\{\alpha_0,\alpha_1,\alpha_2,...\}$ such that
 $$
 \alpha_k\in \ZZ/d_1...d_{k+1}\ZZ\quad \text{and}
 \quad \alpha_k|\text{mod} \,d_1...d_k=\alpha_{k-1}.
 $$
 In particular, if $p$ is a prime number and
 $d_k=p$ we get the ring of $p$--adic integers, $\II=\ZZ_p$.
 
In \cite{PVY1} we built a certain machine that allows
to construct almost periodic Jacobi matrices with singularly
continuous spectrum such that $I_J=\II$.

The key element of the construction is the following 

\begin{theorem} 

Let $\tilde J$ be a Jacobi matrix
with the spectrum on $[-1,1]$. Then the following
Renormalization Equation has a solution 
$J=J(\epsilon,\tilde J)=J(\epsilon,\tilde J; T)$ with
the spectrum on $T^{-1}([-1,1])$:
\begin{equation}\label{3}
V^*_{\epsilon}(z-J)^{-1}V_
{\epsilon}=
(T(z)-\tilde J)^{-1}T'(z)/d,
\end{equation}
where $V_\epsilon|k\rangle=|\epsilon+dk\rangle$, 
$0\le\epsilon\le d-1$.
Moreover, 
 if
$\min_i |t_i|\ge 10$ then 
$$
||J(\epsilon,\tilde J_1)-J(\epsilon,\tilde J_2)||\le
c||\tilde J_1-\tilde J_2||.
$$
with an absolute constant $c<1$ (does not depend of $T$ either ).
\end{theorem}

Let us point out the following two properties
of the function $J(\epsilon,\tilde J; T)$. First,  due to the
commutant relation  $V_\epsilon S=S^d V_\epsilon$ one gets
${J(\epsilon,S^{-m}\tilde J S^m)=S^{-d m}J(\epsilon,\tilde J)S^{d m}}$.
Second, the chain rule holds
$$
J(\epsilon_0, J(\epsilon_1,\tilde J; T_2); T_1)=
J(\epsilon_0+\epsilon_1 d_1,\tilde J; T_2\circ T_1),
$$
where $d_i=\deg T_i$, $0\le\epsilon_i\le d_{i+1}$.

Next steps are quite simple.
For given $d_1, d_2...$, let us chose polynomials
$T_1, T_2...$, $\deg T_k=d_k$ with sufficiently large
critical values. For a fixed sequence $\epsilon_0, \epsilon_1...$,
$0\le\epsilon_k\le d_{k+1}$, 
define
$J_m=J(\epsilon_0+\epsilon_1 d_1+...+\epsilon_{m-1} d_1...d_{m-1},
\tilde J; T_m\circ ...\circ T_2\circ T_1)$. Then $J=\lim_{m\to\infty}
J_m$ exists and does not depend of $\tilde J$. Moreover,
$$
\forall j,\,\,||J-S^{-d_1...d_l j}JS^{d_1...d_l j}||\le A c^l,
\ A>0.
$$
That is $\rho_J$  defines on $\ZZ$ the standard 
$p$--adic topology  in this
case.

Notice that for the case $T_1=T_2=...=T_m=:T$, $T=f^n$, $d_1=...=N^n=:d$ we just get
\begin{equation}
\label{AlmostPeriod}
\forall j,\,\,||J-S^{-d^l j}JS^{d^l j}||\le A c^l,
\ A>0\,.
\end{equation}
This proves that $J$ is a limit periodic matrix (so, in particular) it is almost periodic.
This provides the bridge between  Lipschitz or (better) contractivity property of our noncommutative PFR operator and the question of almost periodicity of a wide class of Jacobi matrix generated by hyperbolic dynamical systems.


\end{document}